# TORSION POINTS OF ABELIAN VARIETIES IN ABELIAN EXTENSIONS

WOLFGANG M. RUPPERT

ABSTRACT. We show that if $A$ is an abelian variety defined over a number field $K$ then $A_{\mathrm{tors}}(K^{\mathrm{ab}})$ is finite iff $A$ has no abelian subvariety with complex multiplication over $K$. We apply this to give another proof for Ribet's result that $A_{\mathrm{tors}}(K^{\mathrm{cycl}})$ is finite.

## 1. INTRODUCTION

For a field $K$ let $\overline{K}$ denote the algebraic closure of $K$, $K^{\mathrm{ab}}$ the maximal abelian extension of $K$ and $K^{\mathrm{cycl}}$ the field obtained by adjoining all roots of unity to $K$. Then $K \subseteq K^{\mathrm{cycl}} \subseteq K^{\mathrm{ab}} \subseteq \overline{K}$.

Let $A$ be an abelian variety defined over a number field $K$ and let $A_{\mathrm{tors}}$ denote the torsion subgroup of $A$. The Mordell-Weil theorem shows that $A_{\mathrm{tors}}(K)$ is finite. Ribet [R] has shown that $A_{\mathrm{tors}}(K^{\mathrm{cycl}})$ is finite. Our aim is to prove the following theorem:

**Theorem.** *Let $A$ be an abelian variety defined over a number field $K$ such that $A$ is $K$-simple. Then $A_{\mathrm{tors}}(K^{\mathrm{ab}})$ is infinite if and only if $A$ has complex multiplication over $K$. In this case $A_{\mathrm{tors}}(\overline{K}) = A_{\mathrm{tors}}(K^{\mathrm{ab}})$.*

We say that an abelian variety $A$ has complex multiplication over $K$ if $\mathrm{End}_K(A) \otimes_{\mathbf{Z}} \mathbf{Q}$ is a number field of degree $2 \dim A$ over $\mathbf{Q}$. An easy consequence of the theorem is the following corollary:

**Corollary 1.** *Let $A$ be an abelian variety defined over a number field $K$. Then $A_{\mathrm{tors}}(K^{\mathrm{ab}})$ is finite if and only if $A$ has no abelian subvariety with complex multiplication over $K$.*

From the theorem we will also deduce another proof of Ribet's result:

**Corollary 2.** *Let $A$ be an abelian variety defined over a number field $K$. Then $A_{\mathrm{tors}}(K^{\mathrm{cycl}})$ is finite.*

The proof of the theorem makes essential use Faltings' finiteness theorems for abelian varieties over number fields.

## 2. PREPARATIONS

Our first lemma is purely algebraic.

**Lemma 1.** *Let $V$ be a finite dimensional vector space over a field $k$. Let $R, S \subseteq \mathrm{End}_k(V)$ be $k$-subalgebras. Let $U \subseteq V$ be a $k$-subspace such that $U \neq 0$. We make the following assumptions:*

---

- $R$ *is a semisimple $k$-algebra.*
- $V$ *is a free $R$-module.*
- $V$ *is a semisimple $S$-module.*
- $R = \mathrm{End}_S(V)$, *i.e. $R$ is the commutant of $S$ in $\mathrm{End}_k(V)$.*
- $RU \subseteq U$ *and* $SU \subseteq U$.
- *The image of the (natural) ring homomorphism $S \to \mathrm{End}_k(U)$ is commutative.*

*Then there is a ring homomorphism $R \to \overline{k}$ (where $\overline{k}$ is the algebraic closure of $k$) and* $\dim_k R = \dim_k V$.

*Proof.* We first prove the lemma in the case that $R$ is a simple $k$-algebra.

As $V$ is a semisimple $S$-module there is a $k$-subspace $W \subseteq V$ such that $SW \subseteq W$ and $V = U \oplus W$. Define $\alpha \in \mathrm{End}_k(V)$ by $\alpha|_U = 0_U$ and $\alpha|_W = 1_W$. Then $\alpha$ commutes with all elements of $S$ and therefore $\alpha \in R$. As $U$ is $R$-invariant we get by restriction a ring homomorphism $R \xrightarrow{\varphi} \mathrm{End}_k(U)$. As $R$ is simple and $\varphi(\alpha) = 0$ we get $\alpha = 0$ and therefore $W = 0$ such that $U = V$. This shows that $S$ itself is commutative. As $R = \mathrm{End}_S(V)$ we get $S \subseteq R$. As $R$ is simple, $S$ is the commutant of $R$ in $\mathrm{End}_k(V)$ (by the density theorem [B, p.39]) and the formula in [B, Théorème 2, p.112] gives then

$$\dim_k R \cdot \dim_k S = (\dim_k V)^2.$$

The fact that $V$ is a free $R$-module gives $\dim_k R \leq \dim_k V$ and $S \subseteq R$ gives $\dim_k S \leq \dim_k R$. Therefore

$$\dim_k R \cdot \dim_k S \leq (\dim_k R)^2 \leq (\dim_k V)^2 = \dim_k R \cdot \dim_k S.$$

This shows $\dim_k R = \dim_k S = \dim_k V$ and $S = R$. Then $R$ is commutative and therefore a field. As $R \subseteq \mathrm{End}_k(V)$ clearly $R$ is algebraic over $k$ such that we have an embedding $R \to \overline{k}$. This proves the lemma in case $R$ is a simple $k$-algebra.

Now we consider the general case.

As $R$ is a semisimple $k$-algebra there are idempotents $\varepsilon_1, \ldots, \varepsilon_r \in R$ such that $\varepsilon_i^2 = \varepsilon_i$, $\varepsilon_i \varepsilon_j = 0$ for $i \neq j$, $1 = \varepsilon_1 + \cdots + \varepsilon_r$ and $R_i = \varepsilon_i R = R\varepsilon_i$ is a simple $k$-algebra. We have $R = R_1 \oplus \cdots \oplus R_r$. $\varepsilon_i$ is the unit element in $R_i$. If we write $V_i = \varepsilon_i V$ we get the decomposition $V = V_1 \oplus \cdots \oplus V_r$. $V_i$ is a $R_i$-module. We also have $\varepsilon_i|_{V_i} = 1_{V_i}$. If $v = v_1 + \cdots + v_r \in V$ with $v_i \in V_i$ then $v_i = \varepsilon_i v$. This implies

$$V_i = \{v \in V : \varepsilon_1 v = \cdots = \varepsilon_{i-1} v = \varepsilon_{i+1} v = \cdots = \varepsilon_r v = 0\}.$$

By assumption $V$ is a free $R$-module: $V \simeq R \oplus \cdots \oplus R = R^\ell$. Then

$$V \simeq (R_1 \oplus \cdots \oplus R_r) \oplus \cdots \oplus (R_1 \oplus \cdots \oplus R_r)$$

and therefore $V_i \simeq R_i^\ell$ such that $V_i$ is also a free $R_i$-module and

$$\frac{\dim_k V_i}{\dim_k R_i} = \ell = \frac{\dim_k V}{\dim_k R}.$$

As $S$ commutes with $R$ we see by the above expression for $V_i$ that $SV_i \subseteq V_i$. Let $S_i$ be the image of $S \to \mathrm{End}_k(V_i)$. Then $S = S_1 \oplus \cdots \oplus S_r$ and $S_i = \varepsilon_i S$. We have $R_i, S_i \subseteq \mathrm{End}_k(V_i)$. It is easy to see that

$$R_i = \mathrm{End}_{S_i}(V_i).$$

As $V$ is a semisimple $S$-module, $V_i$ is a semisimple $S$-module and therefore a semisimple $S_i$-module. Define $U_i = \varepsilon_i U$. Then $U_i \subseteq V_i$ satisfies $R_i U_i \subseteq U_i$ and $S_i U_i \subseteq U_i$. It is clear that the image of the induced ring homomorphism $S_i \to \mathrm{End}_k(U_i)$ is also commutative. As $U = U_1 \oplus \cdots \oplus U_r$ and $U \neq 0$ there is an index $i$ such that $U_i \neq 0$. Now we can apply the first part of the proof and get $\dim_k R_i = \dim_k V_i$ and a ring homomorphism $R_i \to \overline{k}$ which gives by the above formulas $\dim_k R = \dim_k V$ and a ring homomorphism $R \to R_i \to \overline{k}$. ∎

We prove the following lemma for lack of a reference.

**Lemma 2.** *Let $A$ be an abelian variety of dimension $n$ defined over $\mathbf{C}$. Let $\mathcal{O} \subseteq \mathrm{End}(A)$ be a ring of endomorphisms of rank $d$ over $\mathbf{Z}$ such that $D = \mathcal{O} \otimes_{\mathbf{Z}} \mathbf{Q}$ is a division algebra over $\mathbf{Q}$. Then we have:*

1. *$A[p]$ is a free $\mathcal{O} \otimes_{\mathbf{Z}} \mathbf{Z}/(p)$-module of rank $\frac{2n}{d}$ if $p$ is sufficiently large.*
2. *$V_p(A)$ is a free $D \otimes_{\mathbf{Q}} \mathbf{Q}_p$-module of rank $\frac{2n}{d}$ for all $p$.*

*Proof.* There is a lattice $\Lambda \subseteq \mathbf{C}^n$ such that analytically $A \simeq \mathbf{C}^n / \Lambda$. The $p^\ell$-torsion points are then $A[p^\ell] = \frac{1}{p^\ell}\Lambda/\Lambda$ and $T_p(A) = \lim_{\leftarrow} \frac{1}{p^\ell}\Lambda/\Lambda$ with the transition maps $\frac{1}{p^\ell}\Lambda/\Lambda \xleftarrow{\cdot p} \frac{1}{p^{\ell+1}}\Lambda/\Lambda$. Each $\alpha \in \mathcal{O}$ is given by a matrix $M(\alpha) \in M_n(\mathbf{C})$ such that $M(\alpha)\Lambda \subseteq \Lambda$. This gives $\Lambda$ the structure of an $\mathcal{O}$-module. Therefore $\Lambda \otimes_{\mathbf{Z}} \mathbf{Q}$ is a $D$-vector space. By comparing dimensions over $\mathbf{Q}$ we see that $\Lambda \otimes_{\mathbf{Z}} \mathbf{Q}$ has dimension $\frac{2n}{d}$ over $D$. In particular $d|2n$. Let $e_1, \ldots, e_r \in \Lambda$ (with $r = \frac{2n}{d}$) be a $D$-basis of $\Lambda \otimes_{\mathbf{Z}} \mathbf{Q}$. Let $\alpha_1, \ldots, \alpha_d$ be a basis of $\mathcal{O}$ over $\mathbf{Z}$. Then $\alpha_i e_j$, $1 \leq i \leq d, 1 \leq j \leq r$ form a $\mathbf{Q}$-basis of $\Lambda \otimes \mathbf{Q}$. This implies that the $\mathbf{Z}$-module generated by $\alpha_i e_j$, $1 \leq i \leq d, 1 \leq j \leq r$, has finite index $N$ in $\Lambda$. The vectors $\alpha_i e_j \in \Lambda \otimes \mathbf{Q}$ are linearly independent over $\mathbf{Q}$.

1. We look at $A[p] \simeq \frac{1}{p}\Lambda/\Lambda$. This is a $\mathcal{O} \otimes \mathbf{Z}/(p)$-module. Let $f_j$ the the image of $\frac{1}{p}e_j$ in $A[p] = \frac{1}{p}\Lambda/\Lambda$.

*Claim:* $f_1, \ldots, f_r$ are a basis of $A[p]$ over $\mathcal{O} \otimes_{\mathbf{Z}} \mathbf{Z}/(p)$ if $p$ is prime to $N$.

Suppose that we have $\beta_j \in \mathcal{O}$ such that $\sum_{j=1}^r \beta_j f_j = 0$ in $A[p]$. We write $\beta_j = \sum_i m_{ij}\alpha_i$. Then we get

$$\sum_{i,j} m_{ij}\alpha_i \frac{1}{p}e_j \in \Lambda.$$

Every element of $N\Lambda$ is a linear combination of $\alpha_i e_j$ so that we find $n_{ij} \in \mathbf{Z}$ with

$$N \sum_{i,j} m_{ij}\alpha_i \frac{1}{p}e_j = \sum_{i,j} n_{ij}\alpha_i e_j$$

which implies $N m_{ij} = p n_{ij}$. As $p$ is by assumption prime to $N$ we can write $m_{ij} = p\tilde{m}_{ij}$ with $\tilde{m}_{ij} \in \mathbf{Z}$ and therefore

$$\beta_j = \sum_i m_{ij}\alpha_i = p \sum_i \tilde{m}_{ij}\alpha_i \in p\mathcal{O}$$

so that the image of $\beta_j$ in $\mathcal{O} \otimes \mathbf{Z}/(p)$ is 0. This proves

$$(\mathcal{O} \otimes \mathbf{Z}/(p))f_1 \oplus \cdots \oplus (\mathcal{O} \otimes \mathbf{Z}/(p))f_r \subseteq A[p].$$

Comparing dimensions over $\mathbf{Z}/(p)$ shows that we have equality which proves the claim and the first part of the lemma.

2. Now we investigate $T_p(A)$. Define

$$\tilde{e}_j = (\frac{1}{p}e_j, \frac{1}{p^2}e_j, \frac{1}{p^3}e_j, \dots) \in T_p(A) \subseteq V_p(A).$$

*Claim:* $\tilde{e}_1, \ldots, \tilde{e}_r$ are a basis of $V_p(A)$ over $D \otimes_{\mathbf{Q}} \mathbf{Q}_p$.
Suppose that we have $\tilde{\beta}_j \in D \otimes \mathbf{Q}_p$ such that

$$\tilde{\beta}_1 \tilde{e}_1 + \cdots + \tilde{\beta}_r \tilde{e}_r = 0$$

in $V_p(A)$. Then there are $\tilde{m}_{ij} \in \mathbf{Q}_p$ such that

$$\tilde{\beta}_j = \sum_i \tilde{m}_{ij} \alpha_i.$$

By multiplication with a $p$-power we can achieve that all $\tilde{m}_{ij} \in \mathbf{Z}_p$ and that not all $\tilde{m}_{ij}$ are divisible by $p$. We have now

$$\sum_{i,j} \tilde{m}_{ij} \alpha_i \tilde{e}_j = 0.$$

Take $\ell \in \mathbf{N}$ with $p^\ell > N$ and choose $m_{ij} \in \mathbf{Z}$ with $m_{ij} \equiv \tilde{m}_{ij} \bmod p^\ell$. Then

$$\sum_{i,j} m_{ij} \alpha_i \frac{1}{p^\ell} e_j \in \Lambda.$$

Therefore we find $n_{ij} \in \mathbf{Z}$ with

$$N \sum_{i,j} m_{ij} \alpha_i \frac{1}{p^\ell} e_j = \sum_{i,j} n_{ij} \alpha_i e_j.$$

This implies $N m_{ij} = p^\ell n_{ij}$ and with $p^\ell > N$ we get $m_{ij} \equiv 0 \bmod p$, contradicting our assumption. Therefore

$$(D \otimes \mathbf{Q}_p) \tilde{e}_1 \oplus \cdots \oplus (D \otimes \mathbf{Q}_p) \tilde{e}_r \subseteq V_p(A).$$

Comparing dimensions over $\mathbf{Q}_p$ gives equality and the claim follows. This proves the second part of the lemma. ∎

**Lemma 3.** *Let $D$ be a noncommutative division algebra of finite dimension over $\mathbf{Q}$ and $\mathcal{O}$ an order in $D$. Then:*

1. *There is no ring homomorphism $D \otimes_{\mathbf{Q}} \mathbf{Q}_p \to k$ where $k$ is a field ($p$ is arbitrary).*
2. *If $p$ is sufficiently large there is no ring homomorphism $\mathcal{O} \otimes_{\mathbf{Z}} \mathbf{Z}/(p) \to k$ where $k$ is a field.*

*Proof.* 1. A ring homomorphism $D \otimes_{\mathbf{Q}} \mathbf{Q}_p \to k$ would give a homomorphism $D \to D \otimes_{\mathbf{Q}} \mathbf{Q}_p \to k$ and as $D$ is a division algebra an embedding $D \hookrightarrow k$, which contracts the assumption that $D$ is noncommutative.

2. Let $\mathfrak{a} \subseteq \mathcal{O}$ be the ideal generated by all elements of the form $xy - yx$, $x, y \in \mathcal{O}$. As $\mathcal{O}$ is noncommutative we have $\mathfrak{a} \neq 0$ and $\mathfrak{a}$ has finite index in $\mathcal{O}$, i.e. there is a $N \in \mathbf{Z}$, $N \geq 1$ such that $N\mathcal{O} \subseteq \mathfrak{a}$. Suppose that we have a ring homomorphism $\mathcal{O} \otimes \mathbf{Z}/(p) \to k$ where $k$ is a field. Then $k$ has characteristic $p$. Let $\varphi : \mathcal{O} \to \mathcal{O} \otimes \mathbf{Z}/(p) \to k$. Then $\varphi(\mathfrak{a}) = 0$ and as $N \cdot 1_{\mathcal{O}} \in \mathfrak{a}$ we get

$$0 = \varphi(N \cdot 1_{\mathcal{O}}) = N \cdot 1_k$$

so that $p | N$. This shows that for all $p$ with $p > N$ the claim is true. ∎

## 3. Proof of the Theorem

Let $A$ be an abelian variety defined over a number field $K$ which is $K$-simple, i.e. $\operatorname{End}_K(A) \otimes_{\mathbf{Z}} \mathbf{Q}$ is a finite dimensional division algebra over $\mathbf{Q}$. Assume first that $A_{\operatorname{tors}}(K^{\operatorname{ab}})$ is infinite. There are two possible cases:

- There are infinitely many $p$ such that $A[p](K^{\operatorname{ab}}) \neq 0$.
- There is a $p$ such that $\cup_{\ell \geq 1} A[p^\ell](K^{\operatorname{ab}})$ is infinite.

We consider the cases separately and deduce in each case that $A$ has complex multiplication over $K$.

**Case I:** We assume that there are infinitely many $p$ with $A[p](K^{\operatorname{ab}}) \neq 0$. Write $V = A[p]$ and $k = \mathbf{Z}/(p)$. Then $V$ is a $k$-vector space of dimension $2n$. Let $R = \operatorname{End}_K(A) \otimes_{\mathbf{Z}} \mathbf{Z}/(p)$. Then $\dim_k R = d$. We take $p$ large enough such that $R$ can be considered as a $k$-subalgebra of $\operatorname{End}_k(V)$, that $V$ is a free $R$-module of rank $\frac{2n}{d}$ by Lemma 2 and furthermore that $R$ is a semisimple $k$-algebra. Let $G$ be the image of $G_K \to \operatorname{Aut}(A[p])$ and write $S = k[G] \subseteq \operatorname{End}_k(V)$. Taking again $p$ large enough we know by [F, Remarks at the beginning of the proof, p.211] that $V$ is a semisimple $S$-module and $R = \operatorname{End}_S(V)$. Define $U = A[p](K^{\operatorname{ab}})$. Then $U$ is $R$- and $S$-invariant and the image of $S \to \operatorname{End}_k(U)$ is commutative. By our assumption there are infinitely many (large enough in the above sense) $p$ with $U \neq 0$. By Lemma 1 we get $\dim_k R = \dim_k V$, i.e. $d = 2n$ and a ring homomorphism $\operatorname{End}_K(A) \otimes \mathbf{F}_p \to \overline{\mathbf{F}}_p$ for infinitely many $p$. By Lemma 3 this implies that $\operatorname{End}_K(A)$ is commutative and therefore a field. This means that $A$ has complex multiplication over $K$.

**Case II:** We assume that $\cup_{\ell \geq 1} A[p^\ell](K^{\operatorname{ab}})$ is infinite. Write $k = \mathbf{Q}_p$ and $V = V_p(A) = T_p(A) \otimes_{\mathbf{Z}_p} \mathbf{Q}_p$. Define $R = \operatorname{End}_K(A) \otimes_{\mathbf{Z}} \mathbf{Q}_p$ and consider it as $k$-subalgebra of $\operatorname{End}_k(V)$. By Lemma 2 $V$ is a free $R$-module of rank $\frac{2n}{d}$. Let $G$ be the image of $G_K$ in $\operatorname{Aut}(V)$ and $S = k[G] \subseteq \operatorname{End}_k(V)$. By [F, Theorem 1, p.211] we know that $V$ is a semisimple $S$-module and $R$ is the commutant of $S$ in $\operatorname{End}_k(V)$. $T_p(A)$ consists of sequences $(P_\ell)$ such that $P_\ell \in A[p^\ell]$ and $p \cdot P_{\ell+1} = P_\ell$. Define

$$U' = \{(P_\ell)_{\ell \geq 1} \in T_p(A) : K(P_\ell) \subseteq K^{\operatorname{ab}} \text{ for all } \ell \geq 1\}$$

and $U = \mathbf{Q}_p U'$. Then $U$ is a $\mathbf{Q}_p$-vector space and $RU \subseteq U$, $SU \subseteq U$ and the image of $S \to \operatorname{End}_{\mathbf{Q}_p}(U)$ is abelian. It is easy to see that our assumption implies that $U \neq 0$. Lemma 1 gives now $\dim_k R = \dim_k V$, i.e. $d = 2n$, and a ring homomorphism $\operatorname{End}_K(A) \otimes \mathbf{Q}_p \to \overline{\mathbf{Q}_p}$. By Lemma 3 $\operatorname{End}_K(A) \otimes \mathbf{Q}$ is a field. This means that $A$ has complex multiplication over $K$.

Suppose now that $A$ has complex multiplication over $K$, i.e. $F = \operatorname{End}_K(A) \otimes_{\mathbf{Z}} \mathbf{Q}$ is a number field of degree $2 \dim A$ over $\mathbf{Q}$. Let $p$ be any prime. By Lemma 2 $V_p(A)$ is isomorphic to $F \otimes_{\mathbf{Q}} \mathbf{Q}_p$, i.e. there is a $v \in V_p(A)$ such that $V_p(A) = (\operatorname{End}_K(A) \otimes \mathbf{Q}_p)v \simeq \operatorname{End}_K(A) \otimes \mathbf{Q}_p$. Let $G$ be the image of $G_K \to \operatorname{Aut}(V_p(A))$. As $G$ is compatible with endomorphisms $G$ is determined by its operation on $v$ so that we get an injection

$$G \hookrightarrow (\operatorname{End}_K(A) \otimes \mathbf{Q}_p)^* \simeq (F \otimes \mathbf{Q}_p)^*$$

which implies that $G$ is abelian. Therefore $K(\cup_{\ell \geq 1} A[p^\ell]) \subseteq K^{\operatorname{ab}}$. As this holds for all primes we get $A_{\operatorname{tors}}(\overline{K}) = A_{\operatorname{tors}}(K^{\operatorname{ab}})$ as claimed in the theorem.

## 4. Proof of Corollary 1

Let $A$ be an abelian variety defined over a number field $K$.

If $A$ has an abelian subvariety $B$ with complex multiplication over $K$ then $A_{\mathrm{tors}}(K^{\mathrm{ab}}) \supseteq B_{\mathrm{tors}}(K^{\mathrm{ab}})$ which is infinite by the theorem.

Suppose now that $A_{\mathrm{tors}}(K^{\mathrm{ab}})$ is infinite. $A$ is $K$-isogenous to a product $A_1 \times \cdots \times A_r$ of abelian varieties which are defined over $K$ and $K$-simple. Then there is index $i$ such that $(A_i)_{\mathrm{tors}}(K^{\mathrm{ab}})$ is infinite. Therefore $A_i$ has complex multiplication over $K$ by the theorem. The image of the map $A_i \to A_1 \times \cdots \times A_r \to A$ is then an abelian subvariety of $A$ which has complex multiplication over $K$. This proves Corollary 1.

## 5. Proof of Corollary 2

Let $A$ be an abelian variety defined over a number field $K$. We want to show that $A_{\mathrm{tors}}(K^{\mathrm{cycl}})$ is finite. As $A$ is isogenous to a product of $K$-simple abelian varieties we can restrict us to the case that $A$ is $K$-simple, i.e. $\mathrm{End}_K(A) \otimes_{\mathbf{Z}} \mathbf{Q}$ is a finite dimensional division algebra over $\mathbf{Q}$. If $A$ has no complex multiplication over $K$ then $A_{\mathrm{tors}}(K^{\mathrm{ab}})$ is finite (by our theorem) and so is $A_{\mathrm{tors}}(K^{\mathrm{cycl}}) \subseteq A_{\mathrm{tors}}(K^{\mathrm{cycl}})$. Therefore it remains to consider the case that $A$ has complex multiplication over $K$. If necessary we can enlarge the field $K$ or change to a $K$-isogenous abelian variety.

As the argument is very explicit for elliptic curves we start with them. For abelian varieties we can argue in a similar way by using a theorem of Shimura.

5.1. **Elliptic curves.** Let $E$ be an elliptic curve defined over a number field $K$ such that $\mathrm{End}_K(E) \supseteq \mathbf{Z}[\sqrt{d}]$ for some $d < 0$. We can enlarge $K$ such that $K$ is Galois over $\mathbf{Q}$ and $E$ is isogenous to $\mathbf{C}/\mathbf{Z}[\sqrt{d}]$ over $K$. Therefore we can assume that $E \simeq \mathbf{C}/\mathbf{Z}[\sqrt{d}]$. Then $j = j(E)$ can be calculated with $q = e^{2\pi i \cdot \sqrt{d}} = e^{-2\pi\sqrt{|d|}}$ and $\sigma_k(n) = \sum_{d|n} d^k$ as

$$j = 1728 \frac{g_2^3}{g_2^3 - 27g_3^2}$$

where

$$g_2 = \frac{4}{3}\pi^4 + 320\pi^4 \sum_{n=1}^{\infty} \sigma_3(n)q^n, \quad g_3 = \frac{8}{27}\pi^6 - \frac{448}{3}\pi^6 \sum_{n=1}^{\infty} \sigma_5(n)q^n.$$

This implies $j \in \mathbf{R}$. We have $\mathbf{Q}(j, \sqrt{d}) \subseteq K$. Let $K_+$ be the real subfield of $K$. Then $j \in K_+$ and $E$ is defined over $K_+$. But as $\sqrt{d} \notin K_+$ we get $\mathrm{End}_{K_+}(E) = \mathbf{Z}$ and therefore $E_{\mathrm{tors}}(K_+^{\mathrm{ab}})$ is finite by our theorem. As $K = K_+(\sqrt{d})$ we have $K^{\mathrm{cycl}} = K_+(\sqrt{d}, \zeta_\ell, \ell \in \mathbf{N}) \subseteq K_+^{\mathrm{ab}}$ (where $\zeta_\ell = e^{2\pi i/\ell}$) which shows that $E_{\mathrm{tors}}(K^{\mathrm{cycl}})$ is finite.

5.2. **Abelian varieties.** Let $A$ be a $K$-simple abelian variety with complex multiplication defined over a number field $K$, in particular $\dim_{\mathbf{Z}} \mathrm{End}_K(A) = 2 \dim A$. By enlarging $K$ and using an isogenous abelian variety we can achieve the following situation according to a theorem of Shimura [L, p.142, Theorem 6.1]: $A$ is defined over a number field $K$ which is Galois over $\mathbf{Q}$; $K_+ = K \cap \mathbf{R}$ has index 2 in $K$; the abelian variety $A$ is already defined over $K_+$ and $\dim_{\mathbf{Z}} \mathrm{End}_{K_+}(A) < 2 \dim A$. This implies by our theorem that $A_{\mathrm{tors}}(K^{\mathrm{ab}}_+)$ is finite. But as before we have $K^{\mathrm{cycl}} \subseteq K^{\mathrm{ab}}_+$ and the finiteness of $A_{\mathrm{tors}}(K^{\mathrm{cycl}})$ follows.

## References

[B]   N. Bourbaki, Algèbre, Chapitre 8, Modules et anneaux semi-simples, Hermann, Paris 1958.

[F]   G. Faltings, Complements to Mordell, in: G. Faltings, G. Wüstholz et al., Rational Points, Vieweg, Braunschweig 1984.

[KL]  N. Katz, S. Lang, Finiteness theorems in geometric class field theory, L'enseignement Mathématique **27** (1981), 285–314.

[L]   S. Lang, Complex Multiplication, Springer-Verlag 1983.

[R]   K. Ribet, Torsion points of abelian varieties in cyclotomic extensions, L'enseignement Mathématique **27** (1981), 315–319. (Appendix to [KL])

[S]   G. Shimura, On the zeta function of an abelian variety with complex multiplication, Ann. Math. **94** (1971), 504–533.

Institut für Experimentelle Mathematik, Universität GH Essen, Ellernstrasse 29, D-45326 Essen, Germany

*E-mail address*: `ruppert@exp-math.uni-essen.de`